# Bi-level Mixed-Integer Nonlinear Optimization for Pelagic Island Microgrid Group Energy Management Considering Uncertainty

Jichen Zhang, *Student Member IEEE*, Xuan Wei, *Student Member IEEE*, and Yinliang Xu, *Senior Member, IEEE*

*Abstract*—To realize the safe, economical and low-carbon operation of the pelagic island microgrid group, this paper develops a bi-level energy management framework in a joint energy-reserve market where the microgrid group (MG) operator and renewable and storage aggregators (RSA) are independent stakeholders with their own interests. In the upper level, MG operator determines the optimal transaction prices with aggregators to minimize MG operation cost while ensuring all safety constraints are satisfied under uncertainty. In the lower level, aggregators utilize vessels for batteries swapping and transmission among islands in addition to energy arbitrage by participating in energy and reserve market to maximize their own revenue. An upper bound tightening iterative algorithm is proposed for the formulated problem with nonlinear terms and integer variables in the lower level to improve the efficiency and reduce the gap between upper bound and lower bound compared with existing reformulation and decomposition algorithm. Case studies validate the effectiveness of the proposed approach and demonstrate its advantage of the proposed approach in terms of optimality and computation efficiency, compared with other methods.

## NOMENCLATURE

*Abbreviations*

| | |
|---|---|
| DER | Distributed Energy Resource |
| MG | Microgrid Group |
| RSA | Renewable and storage aggregators |
| DER | Distribution Energy Resource |
| EV | Electric Vehicle |
| ESS | Energy Storage System |
| KKT | Karush-Kuhn-Tucker |
| RBRD | Relaxation-based Bi-level Reformulation and Decomposition |

*Variables*

| | |
|---|---|
| $P_{DG,m,i,t}$ | Power output of diesel generator on island $i$ at time slot $t$ at node $m$ |
| $P_{shed,i,t}$ | Power shedding of island $i$ at time slot $t$ |
| $P_{ess,m,i,t}$ | Charging/discharging power of energy storage system on island $i$ at time slot $t$ at node $m$ |
| $P_{agg,m,i,t}$ | Transaction power between RSA and microgrid on island $i$ at time slot $t$ at node $m$ |
| $\pi_{m,i,t}$ | Transaction price of island $i$ at time slot $t$ at node $m$ |
| $\pi^{(\cdot)}_{r,m,i,t}$ | Upward/downward reserve price of island $i$ at time slot $t$ at node $m$ |
| $R^{(\cdot)}_{DG,m,i,t}$ | Upward/downward reserve provided by diesel generators of island $i$ at time slot $t$ at node $m$ |
| $R^{(\cdot)}_{ess,m,i,t}$ | Upward/downward reserve provided by energy storage system of island $i$ at time slot $t$ at node $m$ |
| $R^{(\cdot)}_{agg,m,i,t}$ | Upward/downward reserve provided by RSA to microgrid on island $i$ at time slot $t$ at node $m$ |
| $V_{m,i,t}$ | Nodal voltage of island $i$ at time slot $t$ at node $m$ |
| $\theta_{m,i,t}$ | Nodal phase angle of island $i$ at time slot $t$ at node $m$ |
| $N_{ess,m,i,t}$ | Number of batteries swapping on island $i$ at time slot $t$ at node $m$ |
| $E^{inuse}_{ess,m,i,t}$ | Energy of batteries in use in the energy storage system at node $m$ on island $i$ at time slot $t$ |
| $E_{ess,m,i,t}$ | Energy in the energy storage system at node $m$ on island $i$ at time slot $t$ |
| $N^{full}_{ess,m,i,t}$ | Number of fully-charged batteries in the energy storage system at node $m$ on island $i$ at time slot $t$ |
| $N^{depleted}_{ess,m,i,t}$ | Number of depleted batteries in the energy storage system at node $m$ on island $i$ at time slot $t$ |
| $N^{inuse}_{ess,m,i,t}$ | Number of batteries in use in the energy storage system at node $m$ on island $i$ at time slot $t$ |

*Parameters*

| | |
|---|---|
| $T$ | Length of time slots |
| $I$ | Number of islands |
| $M$ | Number of nodes |
| $a_{m,i}, b_{m,i}, c_{m,i}$ | Coefficient of quadratic, linear, constant term of diesel generators' cost of island $i$ at node $m$ |
| $a_{EM,m,i}, b_{EM,m,i}, c_{EM,m,i}$ | Coefficient of quadratic, linear, constant term of diesel generators' carbon emission of island $i$ at node $m$ |
| $C^{(\cdot)}_{DG,m,i,t}$ | Upward/downward reserve cost diesel generators of island $i$ at node $m$ |

| Symbol | Description |
|---|---|
| $\Delta_{DG,m,i,\max}$ | Maximum ramping rate of diesel generator of island $i$ at node $m$ |
| $k_{shed,m,i}$ | Coefficient of the penalty of power shedding of island $i$ at node $m$ |
| $P_{DG,\min}$ | Minimum limit of power output of diesel generator |
| $P_{DG,\max}$ | Maximum limit of power output of diesel generator |
| $\pi_{m,i,\min}$ | Minimum transaction price of island $i$ at node $m$ |
| $\pi_{m,i,\max}$ | Maximum transaction price of island $i$ at node $m$ |
| $\pi_{r,m,i,\min}^{(\cdot)}$ | Minimum upward/downward reserve price of island $i$ at node $m$ |
| $\pi_{r,m,i,\max}^{(\cdot)}$ | Maximum upward/downward reserve price of island $i$ at node $m$ |
| $P_{ess,m,i,\min}$ | Minimum limit of power output of energy storage system of island $i$ at node $m$ |
| $P_{ess,m,i,\max}$ | Maximum limit of power output of energy storage system of island $i$ at node $m$ |
| $E_{ess,\min}$ | Minimum limit of energy of energy storage system |
| $E_{ess,\max}$ | Maximum limit of energy of energy storage system |
| $EM_{\max}$ | Maximum carbon emission of diesel generators in pelagic island microgrid groups in one day |
| $P_{wind,i,t}$ | Power output of wind turbines on island $i$ at time slot $t$ |
| $x_{mn,i}$ | Branch admittance between node $m$ and $n$ on island $i$ |
| $r_{mn,i}$ | Branch impedance between node $m$ and $n$ on island $i$ |
| $\overline{V}$ | Upper limit of nodal voltage |
| $\underline{V}$ | Lower limit of nodal voltage |
| $N_{ess,m,i}$ | Total number of batteries in the energy storage system at node $m$ on island $i$ at time slot $t$ |
| $E_{ev,j,\max}$ | Maximum limit of energy storage of electric vessel $j$ |
| $E_{ev,j,\min}$ | Minimum limit of energy storage of electric vessel $j$ |
| $n_{v,m,i,t,\max}$ | Maximum number of vessels serving for batteries transmission at node $m$ on island $i$ at time slot $t$ |
| $z$ | Number of batteries each vessel can carry |
| $e$ | Energy of each battery |
| $\xi_{wind,m,i,t}$ | Random variables of wind power at node $m$ on island $i$ at time slot $t$ |
| $\mu_{wind,m,i,t}$ | Mean value of wind power random variables at node $m$ on island $i$ at time slot $t$ |
| $\sigma_{wind,m,i,t}$ | Variance of wind power random variables at node $m$ on island $i$ at time slot $t$ |
| $\xi_{load,m,i,t}$ | Random variables of load at node $m$ on island $i$ at time slot $t$ |
| $\mu_{load,m,i,t}$ | Mean value of load random variables at node $m$ on island $i$ at time slot $t$ |
| $\sigma_{load,m,i,t}$ | Variance of load random variables at node $m$ on island $i$ at time slot $t$ |
| $\varepsilon$ | Violation probability of chance constraints |
| $\Phi$ | Probability distribution function |

## I. Introduction

For a long time, pelagic islands rely on diesel generators to generate electricity, which causes severe environmental pollution and much carbon emission. It is difficult to meet the increasing load demand when diesel generator is used as the only way for power supply. Submarine power cable construction is challenging and expensive since the pelagic island microgrid groups are typically remote from the mainland. Moreover, the size and the local populations of the pelagic island microgrid groups are usually small [12]-[14]. Combined with advanced technology, pelagic island microgrid groups are proposed to solve the integration of renewable energy resources and conventional generators. Due to the characteristic of inhomogeneous spatial distribution of resources and loads, known as resource islands and load islands, the utilization of energy storage system (ESS) and vessels plays an important role in the microgrid groups to achieve the optimal operation. However, few researches focus on the energy management of the pelagic island microgrid groups considering microgrid operator and other stakeholders such as energy storage operators and renewable energy resource operators. When vessels have their own tasks, service fees paid for the vessels are required to be considered in the management.

Researchers mainly focus on the energy management and control strategies in microgrids [1]. In [2], [3], researchers propose privacy-preserve scheduling approaches applied on the distributed energy resources which are sensitive to the local information. In [4], [5], cooperation operation frameworks among multi-microgrids are introduced to guarantee the fairness of profit obtaining. In [6], a two-phase approach is proposed to solve the direct energy trading among multiple microgrids which are a generalized Nash bargaining (GNB) problem that involves the distribution network's operational constraints. In [7], researchers show a state-dimension linear scalable management framework in the microgrids to deal with the energy cost minimization in the multistage stochastic optimization problem. In [8], [9], distributed algorithms for the operation and control of multi-microgrids are introduced which are scalable in larger power systems. Various coordination strategies have been proposed for the distributed microgrids to realize the economic operation [10]-[11], especially in multiple pelagic island microgrid groups. Vessels have been suggested in the energy transmission between the pelagic island microgrid groups. Researchers in [12]-[13] propose an energy management framework by utilizing electric vessels as the way of energy transmission. Sui et al. [14] introduce a day-ahead energy management strategy considering electric vessels dispatching influenced by the weather contributing to the non-integer-hour terms.

TABLE I
COMPARISONS WITH RELATED STUDIES

| Reference | Privacy | Fairness | Scalability | Profitability | | | Non-convexity in Bi-level problem |
| --- | --- | --- | --- | --- | --- | --- | --- |
| | | | | Energy Storage Arbitrage | Reserve Service | Batteries Transmission | |
| [2] | √ | | | √ | | | |
| [3] | √ | | √ | √ | | | |
| [4] | | Cooperative power dispatch | | | | | |
| [5] | | Cooperative power exchange | | √ | | | |
| [6] | | Fair trading price | | | | | |
| [7] | | | √ | | | | |
| [8] | | | Scalable Multi-Microgrids | | | | |
| [9] | | | Scalable Multi-DGs control | | | | |
| [12]-[14] | | | | √ | | √ | |
| [19], [20] | | | √ | | √ | | |
| [29]-[39] | | | | | | | √ |
| Proposed | Limited | Shared benefits | Scalable Multi-Microgrids | √ | √ | √ | √ |

More appropriate frameworks need to be explored when ESS and vessels participant in the energy management. In this paper, renewable and storage aggregators (RSA) are applied for the stable and economic operation of the microgrid groups, cooperating with the vessels. Li *et al.* [21] introduce an aggregation service including wind power, photovoltaic, energy storage and residential loads and model predictive control is applied to realize the optimal DER operations. Duan *et al.* [22] introduce an aggregation of ESS and electric vessels charging stations to bid in both energy and reserve market with conditional value at risk formulated to hedge the risk of randomness. To reach the goal of 'carbon peak, carbon neutrality', carbon emissions need to be reduced in the proposed pelagic island microgrid groups. Wang *et al.* [18] introduce a two-stage scheduling framework for both electricity and carbon markets to mitigate the carbon emissions. Moreover, the interests of other stakeholders are not included in [12]-[14] in the pelagic island groups considering carbon emission and batteries swapping by vessels instead of charging or discharging.

Through further investigation, when considering about the pelagic island, the uncertainty of the renewable energy resources such as wind turbines equipped on the pelagic islands often affects the management of the microgrid, so as the loads. Although previous researches have constructed frameworks for the pelagic islands microgrid groups in the day-ahead scheduling, the uncertainty is not considered exactly in [12]-[14]. To handle with uncertainties of renewable energy resources and loads, reserve services are introduced in most works modeled by chance constrained programming. Wu *et al.* [19] model both reserve service and line flow limits as chance constraints for the day-ahead scheduling with forecast errors in the electricity market. In [20], a chance constrained stochastic optimal dispatch model is formulated to cope with electricity transactions for a wind power integrated energy system and is transformed into a mixed integer second-order cone programming problem. In stochastic programming, the data-based distributions are usually adopted to sample typical scenarios of each type of uncertainty separately, which are combined to form a huge scenario tree. The dispatch strategy is obtained to adapt to these scenarios, which may cause huge computational burdens. In robust programming, the uncertainty set is utilized which is a collection of all possible scenarios. Chance-constrained programming is a tradeoff between the economy and computational burden. Despite that chance constraint programming is frequently used in the reserve market, the application in pelagic island microgrid groups has not been well studied which can ensure the stable power supply under a desired confidence level based on a given distribution, especially considering about the aggregation of renewable energy resources and ESS with vessels serving for batteries swapping and transmission.

Bi-level optimization frequently modeled in the Stackelberg game constructs a framework for such market models by introducing the strategic pricing in the upper-level problem and aggregator in the lower-level problem. In previous electricity market studies, bi-level problems are converted into a single-level problem known as MPEC problems by replacing the lower-level problem with KKT condition. Yi *et al.* [23] introduce a bi-level operation scheme including setting transaction price for both active power and reactive power. Yu *et al.* [25] propose a control strategy based on bi-level optimization for EV charging from spatial and temporal perspective which both benefits the grid and EV drivers, and solve it by KKT reformulation. Nevertheless, KKT condition is feasible only when the lower level problem is convex, leading to evolutionary algorithms widely applied in the bi-level optimization problem, especially consisting of integer variables. Paudel *et al.* [26] introduce iterative algorithms in evolutionary games among buyers and sellers in the peer-to-peer trading in microgrids. Zhao *et al.* [27] apply the hybrid multi-objective differential evolutional algorithm in the formulated bi-level bi-objective optimization problems for the reactive power management of DERs. A few researches have been done to find the global optimum in the bi-level optimization problem with integer lower level variables [3], [29]-[36]. Researchers in [31]-[35] introduce some iterative algorithms by tightening the upper and lower bounds to get the optimal solution. Sun *et al.* [3] propose a decentralized transmission and distribution incentive-compatible pricing mechanism to minimize the uplift payments by applying relaxation based bi-level reformulation and decomposition methods. Shen *et al.* [36] introduce an accelerated algorithm to solve the interactions among the DER aggregator, energy and reserve markets considering security checks in the operation of distribution system. Other pricing strategies have been suggested in [37]-[39] when incorporating with non-convexities. However, few researches focus on the solution to bi-level programming with not only binary or integer variables but also nonlinear terms in the lower level, and duality in the relaxed linear programming cannot be applied. Traditional reformulation and decomposition algorithms always lead to large gaps and slow convergence rates when faced with large-scale variables.

Contrary to the researches mentioned above, this paper proposes a bi-level pelagic island microgrid groups operation framework considering both energy and reserve market by chance constrained programming and RSA, in which vessels serve for batteries transmission rather than power lines to realize the power exchange. Due to the nonconvexities in the lower level of the bi-level problem when batteries can only be transported by vessels based on the number of units, a novel reformulation and decomposition-based algorithm is introduced. The proposed algorithm can give a reasonable solution the bi-level operation with nonlinear terms and integer variables in the lower level. Case studies illustrate the effective of the proposed management framework and algorithm.

Therefore, the main contributions of this paper are summarized below:

(1) In contrast to [12]-[14] only considering a single stakeholder, a bi-level pelagic island microgrid groups energy management framework is proposed including RSA and MG regarded as self-interested participants to maximize the RSA's profits when considering batteries swapping through vessels and minimize the MG overall operation cost, respectively.

(2) Compared with the models in [12]-[14] and methods in [24]-[26] which don't take the uncertainty of renewable energy resources and loads into account, chance constrained programming are introduced to enhance the secure operation of MGs and RSAs in the joint energy and reserve market.

(3) In contrast to existing algorithms which may be difficult to deal with the nonlinear term and integer variables and may lead to nonconvergence results or large gaps, a novel reformulation and decomposition algorithm with upper bound tightening is proposed to solve the model, reducing the optimality gap and improving the computational efficiency.

This paper is organized as follows: Section II introduces the energy management framework for pelagic island microgrid groups. In Section III, the proposed reformulation and decomposition algorithm of the operation framework is presented. Case studies are conducted in Section IV. Finally, conclusions are drawn in Section V.

## II. PELAGIC ISLAND MICROGRID GROUP ENERGY MANAGEMENT FRAMEWORK

Following the framework depicted in the literature [13], a pelagic island microgrid group consists of the resources island and load islands. There are diesel generators and wind turbines equipped in each island microgrid. ESS is regarded as an efficient way to coordinate with the uncertainty of renewable energy resources and loads on the island. To deal with the uneven spatial distribution of energy, vessels sailing among the islands serve for battery transmission. Instead of charging or discharging, battery swapping is applied as the way of energy exchanging with vessels can avoid the waiting time for charging/discharging process and vessels propelled by other power can also serve. In this section, we propose a bi-level energy management framework including microgrid operators and RSA operating as two independent stakeholders for their own interests. In this framework, microgrid operator is responsible for setting transaction price for the aggregators. The operation of microgrid groups will benefit from the framework as it maintains the balance of power supply and consumption on the load islands, reducing the power output of the diesel generators which may lead to high carbon emission. Meanwhile, aggregators benefit from combining the ESSs and wind generators as it realizes the consumption of renewable energy resources and gets profit from joint energy and reserve market transactions. High renewable energy penetration and vital load survival requirements result in a significant reserve need. Due to the physical limits of the generators, the reserve requirements cannot be satisfied so that more reserve service needs to be bought from aggregators in the lower level. In the bi-level framework, chance constrained programming is applied to model the uncertainties of the wind and load prediction in the microgrid groups. Therefore, the bi-level optimization problem including integer variables in the lower level is formulated, which is different from the bi-level continuous optimization models in the context of normal microgrid groups on the mainland, where only continuous variables appear in the lower level and the problems are usually convex.

### A. Lower Level for Aggregator of Renewable Energy Resources and Energy Storage System

The goal of the aggregator is to maximize the profit from the power exchange. The loads on the resource island are less so that the consumption of renewable energy resources is hard to satisfy. Aggregators combine the renewable energy resources and ESS to avoid the wind power curtailment. Vessels will transport goods or people from one island to another in the island groups. In the energy management framework, they can also swap depleted batteries with fully-charged batteries on the resource island and fully-charged batteries with depleted batteries on the load island. Batteries can only be transported based on the number of units. Aggregators can practice arbitrage among the microgrid groups due to the different load demands of different islands by paying service fees to the vessels.

In this level, power exchange profits should be considered for the aggregators. Upper-level MG operators will set transaction price for the lower-level aggregators. Then the operator will decide the exchanging power for profits. The profit of aggregators for all time slots at all nodes in the island groups is outlined as follows:

$$\min_{P_{agg,m,i,t}, N_{ess,m,i,t}, R^{(\cdot)}_{ess,m,i,t}, R^{(\cdot)}_{agg,m,i,t}} f_l =$$

$$\sum_t^T \sum_i^I \sum_m^M (\pi_{agg,m,i,t} P_{agg,m,i,t} + k_{ess} P_{ess,m,i,t}^2 + p_{service,t} N_{ess,m,i,t} \quad (1)$$

$$+ C^+_{ess,r,m,t} R^{+^2}_{ess,m,i,t} + C^-_{ess,r,m,t} R^{-^2}_{ess,m,i,t}$$

$$- \pi^+_{r,m,i,t} R^+_{agg,m,i,t} - \pi^-_{r,m,i,t} R^-_{agg,m,i,t})$$

The first term of (1) denotes the income obtained from the power exchanging. Note that the exchanging power can be positive or negative. The second term denotes the charging/discharging cost of the ESS. The third term shows that the service fees need to be paid to the vessels which is an independent stakeholder of another party. For simplicity, the service fees for each vessel is constant for all time slots and all islands in the microgrid groups. The fourth and fifth term illustrate the cost of the reserve providing to the microgrid including upward reserve and downward reserve. The sixth and seventh term mean the earnings of reserve service from the microgrid groups.

The constraints of the ESS operator are shown as follows:

$$P_{ess,m,i,t} + R^+_{ess,m,i,t} \leq P_{ess,m,i,\max} \quad (2)$$

$$P_{ess,m,i,\min} \leq P_{ess,m,i,t} - R^-_{ess,m,i,t} \quad (3)$$

$$P_{ess,m,i,t} - P_{wind,m,i,t} = P_{agg,m,i,t} \quad (4)$$

$$0 \leq N_{ess,m,i,t} \leq n_{v,m,i,t,\max} z \quad (5)$$

$$N_{ess,m,i} = N^{full}_{ess,m,i,t} + N^{depleted}_{ess,m,i,t} + N^{inuse}_{ess,m,i,t} \quad (6)$$

Eqn.(2) and (3) represent the limits of the power output and the reserve of the ESS. (4) illustrates that the exchanging power consists of wind power generation and ESS power output. (5) denotes that the maximal number of batteries can be transported when the number of vessels at certain time slot is given in the day-ahead scheduling. (6) represents that the number of batteries for ESS on each island is fixed. For the number of batteries, we have following constraints:

$$N^{full}_{ess,m,i,t} = N^{full}_{ess,m,i,t-1} - N_{v,m,i,t} + \left\lfloor \frac{E^{inuse}_{ess,m,i,t-1} + P_{ess,m,i,t}\Delta t}{e} \right\rfloor \quad (7)$$

$$N^{depleted}_{ess,m,i,t} = N^{depleted}_{ess,m,i,t-1} + N_{v,m,i,t} - \left\lfloor \frac{E^{inuse}_{ess,m,i,t-1} + P_{ess,m,i,t}\Delta t}{e} \right\rfloor \quad (8)$$

$$E^{inuse}_{ess,m,i,t} = E^{inuse}_{ess,m,i,t-1} + P_{ess,m,i,t}\Delta t - e \left\lfloor \frac{E^{inuse}_{ess,m,i,t-1} + P_{ess,m,i,t}\Delta t}{e} \right\rfloor \quad (9)$$

$$E_{ess,m,i,t} = E^{inuse}_{ess,m,i,t-1} + N^{full}_{ess,m,i,t} e \quad (10)$$

(7) and (8) show the changes in the number of fully-charged batteries and depleted batteries for ESS on each island. (9) denotes the changes in energy storage on each island. Note that the number of batteries swapping at time slot *t* can be positive or negative when positive means that the batteries are swapped out on the resource islands while negative means ESS on the load islands obtains the batteries from vessels. (10) illustrates the energy left in the ESS.

In the lower level, aggregators handle with the uncertainties of the wind turbines power output modeled by chance constrained programming. For simplicity, we model the distribution of the randomness of wind turbines power output as normal distribution and formulate them as follows:

$$\xi_{wind,m,i,t} \sim N(\mu_{wind,m,i,t}, \sigma^2_{wind,m,i,t}) \quad (11)$$

And the chance constrained is shown below:

$$\Pr\{\sum_m^M (P_{wind,m,i,t} + \xi_{wind,m,i,t} + P_{agg,m,i,t} + R^+_{agg,m,i,t}) \leq \sum_m^M (P_{ess,m,i,t} + R^+_{ess,m,i,t})\} \geq 1 - \varepsilon \quad (12)$$

$$\Pr\{\sum_m^M (P_{wind,m,i,t} + \xi_{wind,m,i,t} + P_{agg,m,i,t} - R^+_{agg,m,i,t}) \geq \sum_m^M (P_{ess,m,i,t} - R^-_{ess,m,i,t})\} \geq 1 - \varepsilon \quad (13)$$

Eqn.(12) and (13) denote that upward and downward reserve provided by the ESS should satisfy the reserve requirements of the wind turbines before selling to the MG.

Note that the decision variables of the lower level contain the numbers of vessels and batteries which are integer variables and the objective function has quadratic terms and bilinear terms.

*B. Upper Level for Microgrid*

The goal of the microgrid groups is to minimize the operation cost including diesel generators, loads shedding and reserve service on the island. The specific mathematical model is given below.

The operational cost of diesel generators can be formulated into a quadratic form as follows:

$$f_{DG,m,i,t} = a_{m,i} P^2_{DG,m,i,t} + b_{m,i} P_{DG,m,i,t} + c_{m,i} \quad (14)$$

Considering about the requirements of reducing the environmental pollution, the carbon emission cost related to the power output of diesel generators can be formulated into a quadratic form as follows:

$$f_{EM,m,i,t} = k_{EM}(a_{EM,m,i} P^2_{DG,m,i,t} + b_{EM,m,i} P_{DG,m,i,t} + c_{EM,m,i}) \quad (15)$$

Thus, the penalty of power shedding is formulated into a quadratic form as follows:

$$f_{shed,m,i,t} = k_{shed,m,i} P^2_{shed,m,i,t} \quad (16)$$

Microgrid exchanges power with RSA to get profit or balance the power supply and consumption. The power exchange payment which is opposite to the lower level can be formulated into a linear form as follows:

$$f_{agg,m,i,t} = -\pi_{m,i,t} P_{agg,m,i,t} \quad (17)$$

To cooperate with the uncertainty of the renewable resources and loads, reserve service needs to be considered. The cost of diesel generator reserve including upward and downward reserve is given as follows:

$$f_{DG,r,m,i,t} = C^+_{DG,r,m,t} R^+_{DG,m,i,t} + C^-_{DG,r,m,t} R^-_{DG,m,i,t} \quad (18)$$

The reserve service buying from the RSA which is opposite to the lower level is shown below:

$$f_{agg,r,m,i,t} = \pi^+_{r,m,i,t} R^+_{agg,m,i,t} + \pi^-_{r,m,i,t} R^-_{agg,m,i,t} \quad (19)$$

Therefore, the whole cost of the microgrid including the terms mentioned above is formulated as follows:

$$\min_{P_{DG,m,i,t}, P_{shed,m,i,t}, \pi_{m,i,t}, \pi^{(\cdot)}_{r,m,i,t}, R^{(\cdot)}_{DG,m,i,t}} f_{up} =$$
$$\sum_{t}^{T} \sum_{i}^{I} \sum_{m}^{M} (f_{DG,m,i,t} + f_{EM,m,i,t} + f_{shed,m,i,t}$$
$$+ f_{agg,m,i,t} + f_{DG,r,m,i,t} + f_{agg,r,m,i,t}) \quad (20)$$

Also, the security of operation solution needs to be guaranteed, so that the following constraints should be considered in the model:

$$P_{DG,m,i,\min} \leq P_{DG,m,i,t} - R^-_{DG,m,i,t} \quad (21)$$

$$P_{DG,m,i,t} + R^+_{DG,m,i,t} \leq P_{DG,m,i,\max} \quad (22)$$

$$-\Delta_{DG,m,i,\max} \leq P_{DG,m,i,t} - P_{DG,m,i,t-1} \leq \Delta_{DG,m,i,\max} \quad (23)$$

$$0 \leq P_{shed,m,i,t} \leq P_{load,m,i,t} \quad (24)$$

$$\sum_{t}^{T} \sum_{i}^{I} \sum_{m}^{M} (a_{EM,m,i} P^2_{DG,m,i,t} + b_{EM,m,i} P_{DG,m,i,t} + c_{EM,m,i}) \leq EM_{\max} \quad (25)$$

$$\pi_{m,i,\min} \leq \pi_{m,i,t} \leq \pi_{m,i,\max} \quad (26)$$

$$\pi^{(\cdot)}_{r,m,i,\min} \leq \pi^{(\cdot)}_{r,m,i,t} \leq \pi^{(\cdot)}_{r,m,i,\max} \quad (27)$$

$$P_{m,i,t} = P_{DG,m,i,t} - P_{agg,m,i,t} - P_{load,m,i,t} + P_{shed,i,t} \quad (28)$$

$$P_{m,i,t} = \sum_{m=1, m\neq n}^{b} (\frac{x_{mn,i}}{r^2_{mn,i} + x^2_{mn,i}} (\theta_{m,i,t} - \theta_{n,i,t})$$
$$+ \frac{r_{mn,i}}{r^2_{mn,i} + x^2_{mn,i}} (V_{m,i,t} - V_{n,i,t})) \quad (29)$$

$$Q_{m,i,t} = \sum_{m=1, m\neq n}^{b} (-\frac{r_{mn,i}}{r^2_{mn,i} + x^2_{mn,i}} (\theta_{m,i,t} - \theta_{n,i,t})$$
$$+ \frac{x_{mn,i}}{r^2_{mn,i} + x^2_{mn,i}} (V_{m,i,t} - V_{n,i,t})) \quad (30)$$

$$\underline{V} \leq V_{m,i,t} \leq \overline{V} \quad (31)$$

Eqn(21)-(22) represent the magnitude constraint and reserve constraint of the diesel generator. Eqn(23) denotes the ramping limits of diesel generators. Eqn(24) denotes the constraint of load shedding. Eqn(25) shows the limits of carbon emissions. Eqn(26) shows the limits of transaction price. Eqn(28) represents the power balance constraint of the island microgrid. Eqn(27) shows the limits of reserve price. Eqn(29)-(31) denotes the network constraints on a certain island.

Chance constrained programming in the upper level is used to handle with the uncertainties of loads. The uncertainties are formulated as follow:

$$\xi_{load,m,i,t} \sim N(\mu_{load,m,i,t}, \sigma^2_{load,m,i,t}) \quad (32)$$

And the chance constrained is shown below:

$$\Pr\{\sum_{m}^{M}(P_{DG,m,i,t} - P_{agg,m,i,t} + R^+_{DG,m,i,t} + R^+_{agg,m,i,t}) \geq$$
$$\sum_{m}^{M}(P_{d,m,i,t} - P_{shed,m,i,t} + \xi_{load,m,i,t})\} \geq 1-\varepsilon \quad (33)$$

$$\Pr\{\sum_{m}^{M}(P_{DG,m,i,t} - P_{agg,m,i,t} - R^-_{DG,m,i,t} - R^-_{agg,m,i,t}) \leq$$
$$\sum_{m}^{M}(P_{d,m,i,t} - P_{shed,m,i,t} + \xi_{load,m,i,t})\} \geq 1-\varepsilon \quad (34)$$

Eqn.(33) and (34) denote that the reserve service provided by the diesel generator and RSA should satisfy the load demand with a prescribed violation probability.

Overall, the formulation in the energy management framework consists of two levels: the upper level which minimizes the cost for the pelagic island microgrid groups and the lower level which maximizes the profits for the RSA. The microgrid groups in the upper level set the transaction price for the lower level and then the aggregator in the lower level can buy energy at a low-price island and sell at a high-price island to gain profits by paying fees to the vessels to transport batteries. Therefore, energy storage systems and renewable energy resources are aggregated to access the market indirectly. As shown in Fig. 1, The aggregators make the most profitable scheduling according to the transaction price, and the microgrid operator determines the optimal pricing and energy management schemes after receiving the scheduling results from the aggregators. KKT-based reformulation and distributed algorithms, which are mainly used for the microgrid groups connected by power lines in a convex optimization problem with continuous decision variables, are not appropriate in the pelagic island microgrid groups operation framework due to the integer variables and nonlinear terms in the lower level. Therefore, a novel reformulation and decomposition based algorithm will be investigated in the next section.

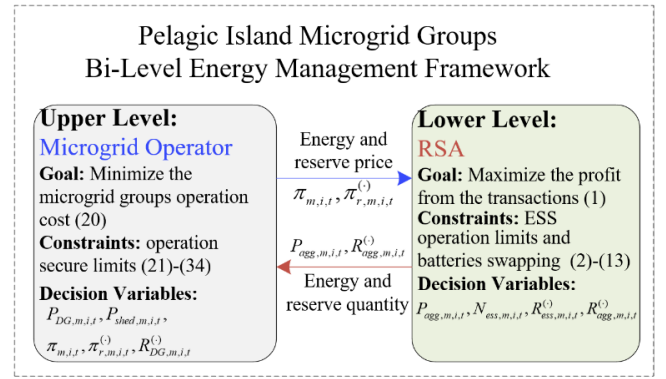

Fig. 1. Bi-level Pelagic Island Microgrid Group Energy Management Framework

## III. BI-LEVEL REFORMULATION AND DECOMPOSITION ALGORITHM WITH NONLINEARITIES

### A. Deterministic Equivalent of Chance Constraint

The chance constraints in the upper and lower level can be converted into the following deterministic constraint:

$$\sum_{m}^{M}(P_{DG,m,i,t} + P_{agg,m,i,t} + R^+_{DG,m,i,t} + R^+_{agg,m,i,t}) \\ - \sum_{m}^{M}(P_{d,m,i,t} - P_{shed,m,i,t}) \geq \Phi(\varepsilon)\sigma_{load,m,i,t} \quad (35)$$

$$\sum_{m}^{M}(P_{DG,m,i,t} + P_{agg,m,i,t} - R^-_{DG,m,i,t} - R^-_{agg,m,i,t}) \\ - \sum_{m}^{M}(P_{d,m,i,t} - P_{shed,m,i,t}) \leq \Phi(\varepsilon)\sigma_{load,m,i,t} \quad (36)$$

$$\sum_{m}^{M}(P_{wind,m,i,t} + P_{agg,m,i,t} + R^+_{agg,m,i,t}) \\ - \sum_{m}^{M}(P_{ess,m,i,t} + R^+_{ess,m,i,t}) \geq \Phi(\varepsilon)\sigma_{wind,m,i,t} \quad (37)$$

$$\sum_{m}^{M}(P_{wind,m,i,t} + P_{agg,m,i,t} - R^-_{agg,m,i,t}) \\ - \sum_{m}^{M}(P_{ess,m,i,t} - R^-_{ess,m,i,t}) \leq \Phi(\varepsilon)\sigma_{wind,m,i,t} \quad (38)$$

### B. Problem Description

The proposed framework in the Section II after deterministic equivalence of chance constraint is essentially a bi-level mixed integer nonlinear programming when there are integer variables and nonlinear term in the lower level. The compact form of the problem can be written as below:

$$\min \sum_{t=0}^{T}\sum_{i=0}^{n} f_{u,i,t}(x_{u,i,t}, x_{l,i,t}, y_{l,i,t}) \\ s.t. Ax_{u,i,t} + Bx_{l,i,t} + Cy_{l,i,t} \leq \mathbf{b}_u \\ x_{l,i,t}, y_{l,i,t} \in \arg\min\{f_{l,i,t}(x_{u,i,t}, x_{l,i,t}, y_{l,i,t}): \\ Dx_{l,i,t} + Ey_{l,i,t} \leq \mathbf{b}_l\}. \quad (39)$$

where $x_{u,i,t}$ represents the continuous variables in the upper level, $x_{l,i,t}$ represents the continuous variables in the lower level and $y_{l,i,t}$ represents the integer variables in the lower level. Note that in the proposed model, variables in the upper level are not contained in the constraints of lower level so that the proposed algorithm can be applied. KKT-based reformulation cannot solve the problem directly.

### C. Reformulation and Decomposition with Nonlinearities

The problem cannot be solved directly by off-the-shelf solvers such as Gurobi due to the constraint(5)-(8), so we need to relax the problem into the following form:

$$\min \sum_{t=0}^{T}\sum_{i=0}^{n} f_{u,i,t}(x_{u,i,t}, x_{l,i,t}, y_{l,i,t}) \\ s.t. Ax_{u,i,t} + Bx_{l,i,t} + Cy_{l,i,t} \leq \mathbf{b}_u \\ f_{l,i,t}(x_{u,i,t}, x_{l,i,t}, y_{l,i,t}) \leq \\ \min\{f_{l,i,t}(x_{u,i,t}, x_{c,l,i,t}, y_{c,l,i,t}): Dx_{c,l,i,t} + Ey_{c,l,i,t} \leq \mathbf{b}_l\}. \quad (40)$$

Because the lower-level program is a minimization problem, constraint in (40) ensures that $(x_{l,i,t}, y_{l,i,t})$ is an optimal solution to the lower-level problem for any given $x_{u,i,t}$. In some simple cases, the equivalent can be handled through enumerating all feasible combinations of integer variables in the lower level.

$$\min \sum_{t=0}^{T}\sum_{i=0}^{n} f_{u,i,t}(x_{u,i,t}, x_{l,i,t}, y_{l,i,t}) \\ s.t. Ax_{u,i,t} + Bx_{l,i,t} + Cy_{l,i,t} \leq \mathbf{b}_u \\ f_{l,i,t}(x_{u,i,t}, x_{l,i,t}, y_{l,i,t}) \leq \\ \min\{f_{l,i,t}(x_{u,i,t}, x_{c,l,i,t}, y_{c0,l,i,t}): Dx_{c,l,i,t} + Ey_{c0,l,i,t} \leq \mathbf{b}_l\} \\ y_{c0,l,i,t} \in Z_0. \quad (41)$$

where $Z_0$ represents the set of all feasible integer variables combinations. Different from [35], we don't need to consider about the violations of the constraints in the lower level for that there are no upper-level variable in the lower-level constraints. However, it is difficult to enumerate all integer variables which may lead to large number of constraints in our proposed framework. It may consume a lot of time so that accelerated decomposition-based iterative method is given. The problem is decomposed into master problem (MP) and subproblem (SP). The master problem in iteration k is illustrated as follows:

$$\min \sum_{t=0}^{T}\sum_{i=0}^{n} f_{u,i,t}(x_{u,i,t}, x_{l,i,t}, y_{l,i,t}) \\ s.t. Ax_{u,i,t} + Bx_{l,i,t} + Cy_{l,i,t} \leq \mathbf{b}_u \\ Dx_{l,i,t} + Ey_{l,i,t} \leq \mathbf{b}_l \\ f_{l,i,t}(x_{u,i,t}, x_{l,i,t}, y_{l,i,t}) \leq \\ \min\{f_{l,i,t}(x_{u,i,t}, x_{c,l,i,t}, y_{c,k,l,i,t}): Dx_{c,l,i,t} + Ey_{c,k,l,i,t} \leq \mathbf{b}_l\} \\ y_{c,k,l,i,t} \in Z_k^l. \quad (42)$$

where $Z_k^l$ denotes the set of integer variable combinations in iteration k. The optimal objective function value gives a lower bound of the problem and the optimal solution of upper level $x^*_{u,i,t}$ is obtained. Then the SP for certain fixed $x^*_{u,i,t}$ is given below.

$$\min \sum_{t=0}^{T}\sum_{i=0}^{n} f_{l,i,t}(x^*_{u,i,t}, x_{l,i,t}, y_{l,i,t}) \\ s.t. Dx_{l,i,t} + Ey_{l,i,t} \leq \mathbf{b}_l. \quad (43)$$

The optimal value $\theta_k(x^*_{u,i,t})$ and solution $(x^*_{l,i,t}, y^*_{l,i,t})$ of lower level in the iteration k are obtained. In the iterative process, another subproblem SP2 is considered to check whether the solution will give a feasible solution to the bi-level problem and the optimal value will give an upper bound for the problem.

$$\min \sum_{t=0}^{T}\sum_{i=0}^{n} f_{u,i,t}(x^*_{u,i,t}, x_{l,i,t}, y_{l,i,t})$$

$$s.t. Ax^*_{u,i,t} + Bx_{l,i,t} + Cy_{l,i,t} \leq \mathbf{b}_u \quad (44)$$

$$Dx_{l,i,t} + Ey_{l,i,t} \leq \mathbf{b}_l$$

$$f_{l,i,t}(x^*_{u,i,t}, x_{l,i,t}, y_{l,i,t}) \leq \theta_k(x^*_{u,i,t}).$$

Similar to the algorithm in the [35], we apply the KKT-based tightening constraints. For the lower-level problem with nonlinear term, the primal and dual feasibility constraints are given below for the continuous lower level variables when the integer variables are fixed.

$$f_{l,i,t}(x_{u,i,t}, x_{l,i,t}, y_{l,i,t}) \leq f_{l,i,t}(x_{u,i,t}, x_{k,l,i,t}, y_{k,l,i,t})$$

$$Dx_{k,l,i,t} + Ey_{k,l,i,t} \leq \mathbf{b}_l$$

$$D^T \pi^k \geq \nabla_{x_{l,i,t}} f_{l,i,t}(x_{u,i,t}, x_{k,l,i,t}, y_{k,l,i,t})^T \quad (45)$$

$$\nabla_{x_{l,i,t}} f_{l,i,t}(x_{u,i,t}, x_{k,l,i,t}, y_{k,l,i,t}) x_{k,l,i,t} =$$

$$-(\pi^k)^T (Dx_{k,l,i,t} + Ey_{l,i,t} - \mathbf{b}_l) + (\pi^k)^T Dx_{k,l,i,t}.$$

Therefore, the Mod-MP problem can be rewritten as below:

$$\min \sum_{t=0}^{T}\sum_{i=0}^{n} f_{u,i,t}(x_{u,i,t}, x_{l,i,t}, y_{l,i,t})$$

$$s.t. Ax_{u,i,t} + Bx_{l,i,t} + Cy_{l,i,t} \leq \mathbf{b}_u$$

$$Dx_{l,i,t} + Ey_{l,i,t} \leq \mathbf{b}_l$$

$$f_{l,i,t}(x_{u,i,t}, x_{l,i,t}, y_{l,i,t}) \leq f_{l,i,t}(x_{u,i,t}, x_{k,l,i,t}, y_{k,l,i,t})$$

$$Dx_{k,l,i,t} + Ey_{k,l,i,t} \leq \mathbf{b}_l \quad (46)$$

$$D^T \pi^k \geq \nabla_{x_{l,i,t}} f_{l,i,t}(x_{u,i,t}, x_{k,l,i,t}, y_{k,l,i,t})^T$$

$$\nabla_{x_{l,i,t}} f_{l,i,t}(x_{u,i,t}, x_{k,l,i,t}, y_{k,l,i,t}) x_{k,l,i,t} =$$

$$-(\pi^k)^T (Dx_{k,l,i,t} + Ey_{l,i,t} - \mathbf{b}_l) + (\pi^k)^T Dx_{k,l,i,t}$$

$$y_{k,l,i,t} \in Z_k^l.$$

### C. Upper Bound Tightening

To accelerate the iterative process, the upper bound will be given by a feasible solution which may be close to optimal solution.

Similar to the algorithm of strongly convex lower-level problem, the proposed algorithm relaxes the integer variables into continuous variables firstly. Then the KKT reformulation can be applied to the bilevel framework.

$$\frac{\partial L}{\partial \tilde{x}_{l,i,t}} = 0 \quad (47)$$

$$0 \leq g(\tilde{x}_{l,i,t}) \perp \mu \geq 0 \quad (48)$$

where $\tilde{x}_{l,i,t}$ represents the relaxed integer variables in the lower level. The relaxed bi-level problem can be converted into a single-level problem by adding constraints (47) and (48) to the upper level while the solution to the relaxed integer variables of the single-level problem may not be integers. Consequently, the single-level problem becomes a mixed integer linear programming problem, and then branch and bound can be applied. The single-level problem will be solved by the iteration of branch and bound so that the feasibility of the results can be guaranteed. When the solution is obtained, a feasible integer variables combination can be added into the integer set. Then, the objective of the relaxation and KKT reformulation can be regarded as the upper bound for the iterative process. Therefore, RBRD algorithm can be accelerated through the upper-bound tightening.

When solving the bi-level problem by the proposed algorithm, global information is needed to carry out the upper-bound tightening and the scheduling is determined by microgrid operators based on the response of the aggregators. Therefore, the aggregators report the operation parameters to the upper level, and then microgrid operators can solve the bi-level optimization problem. According to the assumptions of the proposed management framework, the parameters are accurate and actual, and the transaction energy and reserve service between the upper level and lower level obtained by solving the problem are the same as the result optimized by the aggregators themselves according to the transaction prices.

The complete pseudocode of proposed algorithm is presented in **Algorithm 1**:

---
**Algorithm 1** The Proposed Algorithm
---
1: Initialize $LB = -\infty$, $\xi = 0$, $k = 0$, $Z_0^L = \varnothing$, $UB$ obtained from upper bound tightening.
2: **While** $UB - LB > \xi$ **do**
3:   Solve (46) for $Z_k^L$ to obtain optimal solution $x^*_{u,i,t}$ and set optimal objective function value as $LB$.
4:   Solve (43) for $x^*_{u,i,t}$ to obtain the optimal value $\theta_k(x^*_{u,i,t})$ and solution $(x^*_{l,i,t}, y^*_{l,i,t})$.
5:   Solve (44) for $x^*_{u,i,t}$ and $\theta_k(x^*_{u,i,t})$
6:   **If** (44) is feasible **then**
7:     Obtain the optimal solution $(x_{l,i,t}, y_{l,i,t})$
8:     $UB \leftarrow \min\{UB, \text{ optimal value of } (44)\}$
9:     $y_{l,k} \leftarrow y_{l,i,t}$
10: **else**
11:    $y_{l,k} \leftarrow y^*_{l,i,t}$
12:  $Z_{k+1}^L \leftarrow Z_k^L \cup \{y_{l,k}\}$
13:  $k \leftarrow k+1$, go to step 2
14: Terminate and output the $(x^*_{u,i,t}, x^*_{l,i,t}, y^*_{l,i,t})$ and optimal value of upper level and lower level.
---

## IV. CASE STUDIES

### A. Test System

The performance of the proposed energy management approach for day ahead scheduling is tested on the three-pelagic-island microgrid group. Case studies have been carried out on a computer with a 4-core 3.40 GHz Intel i5-7500 processor and 8 GB of RAM. The proposed approach has been implemented based on Gurobipy in Python 3.7.

As shown in Fig. 2, each island has constructed a microgrid, which consists of diesel generators, wind turbines and ESSs. In addition, there are some vessels responsible for battery transmission to enable energy sharing among different microgrids. Suppose that island 1 is a load island equipped with a 4-bus microgrid, and the other two are resource islands with a simplified microgrid topology where all the resources are connected in parallel to the common point. The day-ahead optimization cycle is 1 hour and the whole scheduling period is 24 hours. The violation probability is set to 5% and the convergence tolerance is set to 0.01.

The profiles of load demand and wind generation are shown in Fig. 3. The number of available vessels for batteries transmission at each time slot is given in the Fig. 4. The transmission from each resource island to the load island by vessels takes two hours. The diesel generator cost coefficients are \$4.05/MW$^2$, \$38.64/MW, \$12.45, and the associated carbon emission coefficients are \$0.1/MW$^2$, \$0.4/MW. The carbon emission cost coefficient is \$15/tonCO$_2$, the load cost coefficient is \$250/MW, the ESS operation cost coefficient is \$24.45/MWh$^2$, the diesel generator reserve cost is \$91.5/MW and the ESS reserve cost is \$31.5/MW$^2$. The service fee of each battery transmission is \$10 per time and each vessel can carry one battery with 0.15MWh. The standard deviations of wind power output and load demand prediction errors are 15% and 10%, respectively. The upper and lower limits of diesel generators (or ESSs) power output are 2.25MW and 0MW (or 1.875MW and 0MW), respectively. The ramping limit of diesel generators is 0.75MW/h. The maximum and minimum energy of ESS are 1.125MWh and 10.125MWh.

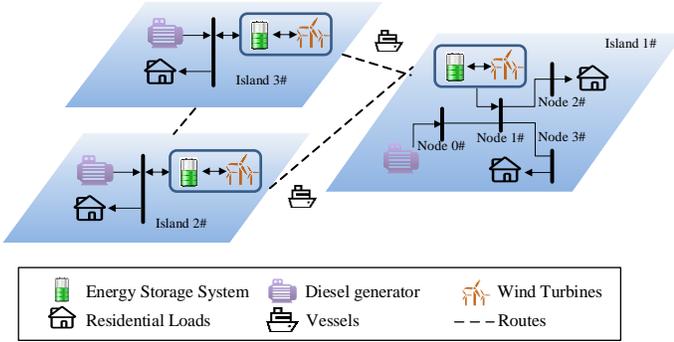

**Fig. 2.** Topology of the three-pelagic-island microgrid group with vessels

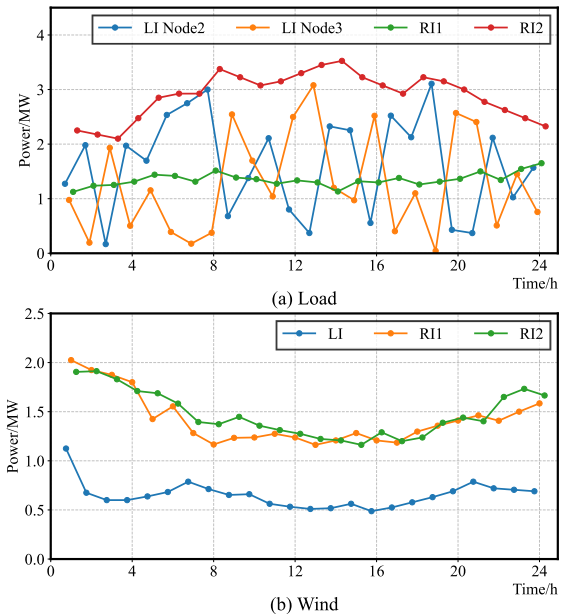

**Fig. 3.** Profiles of predicted load demand (a) and wind power generation (b).

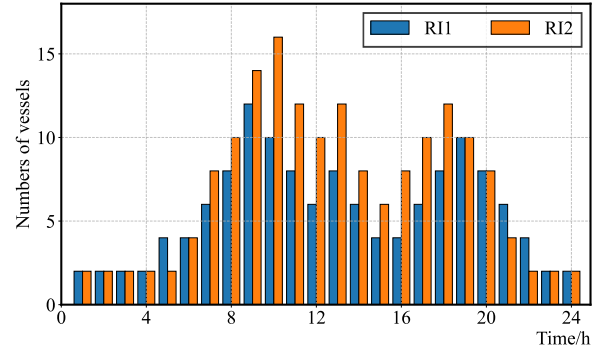

**Fig. 4.** The number of available vessels on three islands.

### B. Simulation Results Analysis

This subsection presents the simulation results of the proposed energy management for the pelagic island microgrids group day-ahead scheduling.

The cost for microgrid group in the upper-level is \$11933.7, and the operation profit for RSA in the lower level is \$5278.6. Fig. 5(a) illustrates the power output of diesel generators on each island in the test system. The diesel generators almost generate no electricity on the resource islands. Compared with the case without ESS and wind turbines equipped on the islands, less diesel generator power output contributes to the lower cost and less environmental pollution.

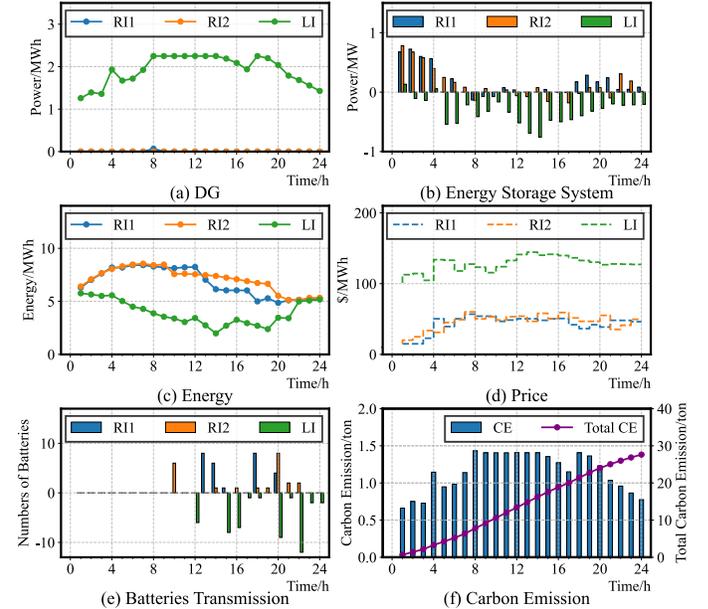

**Fig. 5.** Detailed information of scheduling

Fig. 5(b) shows the batteries power exchange between the microgrid operator and RSA. Fig. 5(c) illustrates the variations of energy in the ESS. The transaction price is shown in Fig. 5(d). Fig. 5(e) shows the number of batteries swapping of each island at each time slot. As it can be seen, the ESS on the resource islands are usually charged while on the load island is usually discharged. The wind turbine generation can be fully consumed by the swapped in empty batteries by the vessels in instead of wind power curtailment due to the fully-charged batteries that need to be swapped out. By paying service fees to vessels, they will serve for the batteries swapping. When the energy of ESS on the load island cannot satisfy the load demand or reach the minimal limitation, the batteries will be swapped from the resource island to the load island and the number is related to the vacancies of the energy on the load island. Aggregators can get profits from arbitrage because of the difference of price between different islands. As is observed, due to the imbalance between the generation and load demand, the price on the load island is always higher than the price on the resource islands. The energy arbitrage can be practiced not only from the temporal perspective but also the spatial perspective. The batteries can store the redundant energy when the load demand is lower than the power output of wind turbines at night, that the price is low, while discharge in the at noon, that the price is high. Moreover, the results show that there is no load shedding in the proposed operation framework. Fig. 5(f) represents the carbon emission at time slot *t* for the whole microgrid groups. The carbon emission limit is set to 30 tons and the proposed approach can satisfy the carbon emission constraint.

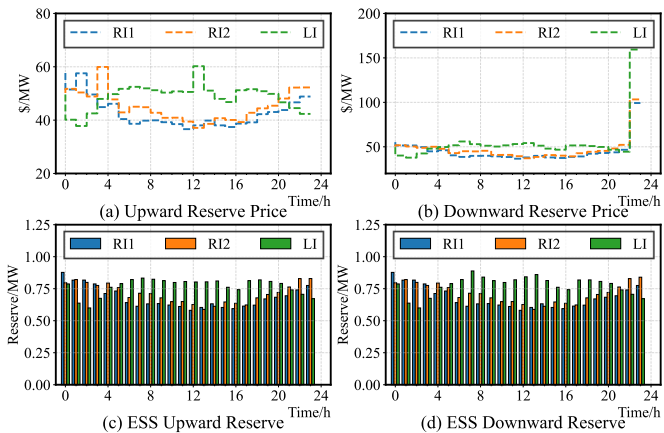

**Fig. 6.** Detailed information of reserve service

Fig. 6 represents the reserve prices which include upward and downward prices provided from the microgrid operator to the aggregators. Since the reserve cost provided by the diesel generators is high, microgrid operator usually tends to buy reserve from the aggregator so the reserve provided by the diesel generators is always zero. The reserve price corresponds with the total reserve requirements including the uncertainty of load and renewable energy resources on the island. On the load island, the uncertainty of the load demand is greater than the uncertainty of renewable energy resources in daytime so that the reserve price on the load island is a little bit higher than the resource island. While the uncertainty of wind turbines on the resource island is greater than the uncertainty of load demand at night, it is obvious that the corresponding reserve price is higher, so as the reserve service provided by ESS. It should be noted that ESS should guarantee that there is enough energy for the operation of the next day which means that the minimal energy of ESS at the last time slot is higher than other time slots so that the downward reserve price at the end of the day for the ESS may be higher than other time slots.

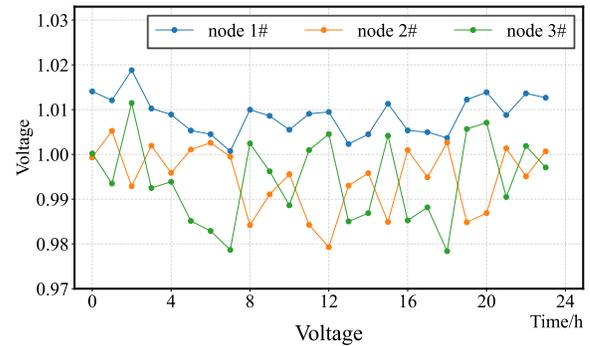

**Fig. 7.** Nodal voltages on load island

Fig. 7 shows that the voltage of load nodes on the load island. Through the control of reactive power in the proposed energy management framework, the nodal voltage is always in the safe range.

*C. Comparison with other methods*

The convergence of the proposed approach is investigated in this subsection. To accelerate the convergence speed, the result of B&B relaxed KKT reformulation is applied as the upper bound of the iterative process.

TABLE II
COMPARISONS AMONG THE THREE ALGORITHMS

| Algorithm | Time (s) | Gap | MG Cost | Iteration Times | Privacy |
|---|---|---|---|---|---|
| EA | 10000 | / | 18432.2 | 95 | √ |
| RBRD | 12000 | 17.67% | 14436.3 | 360 | √ |
| Proposed | 3095 | 0.97% | 11933.7 | 105 | Limited |

As shown in Tab. II, the computation time of the proposed approach is 3095s, which is much less comparing with those of the EA and RBRD algorithm. The iteration time of the proposed algorithm is comparable to that of EA, while much less than RBRD algorithm. Furthermore, the gap between the upper level and lower level in the iterative process of the proposed algorithm is 0.97% when the iteration is terminated which is smaller than the RBRD algorithm, when EA focuses on searching the minimal solution instead of tightening the upper bound and lower bound so that the comparison of gap cannot be conducted. In the proposed algorithm, the tightened upper bound and a feasible integer variables combination of the operation framework derived from the relaxed KKT reformulation and branch and bound have been obtained while in RBRD algorithm, the upper bound is infinite and the primal integer variables combination set is empty so that the computational time and the gap between the upper and lower bound can be significantly reduced. The proposed algorithm also achieves a lower cost $11933.7 of the operation framework, when the RBRD cannot give a convergence result and EA is difficult to find out the minimum in such a high-dimensional variable space. However, to tighten the bound of the iterative process which can accelerate the algorithm, the information of aggregators needs to be reported to microgrid operators so that privacy cannot be preserved while the other two algorithms can protect the private information of each stakeholder. The result demonstrates the effectiveness of the proposed algorithm, showing great edge on the convergence and solution optimality.

To validate the economic performance of reserve service and RSA, the comparison between the framework and other cases is conducted and the result is shown in Tab III. In the table, 'time' refers to the computation time for the day-ahead optimal operation on the platform mentioned in the Part A and the 'total cost' refers to the operation cost including the energy cost and reserve cost under three different strategies.

TABLE III
COMPARISON BETWEEN THE TWO CONDITIONS

|  | Proposed | Without reserve | Without RSA |
| --- | --- | --- | --- |
| Total cost ($) | 11933.74 | 9681.48 | 14280.47 |
| Reliability | 95.4% | 0% | 80.7% |
| Time(s) | 3095 | 670 | 1 |

The first case is that the reserve service for diesel generators and ESS is excluded from the framework and the second case is that there is only diesel generators and wind turbines equipped instead of RSA and reserve service is provided by diesel generators. In comparison, Monte Carlo simulation is introduced to generate massive scenarios when the distribution of uncertainty is given. When uncertainty is not considered and reserve service is not provided, the operation framework is unreliable with the 100% simulation constraint violations occur despite of lower cost. When RSA is not considered, the cost may be high due to the load shedding and the reliability may be worse due to the limitations of reserve service provided by diesel generators. In the proposed operation framework, only 4.6% violations occur with the reserve service dealing with the uncertainty.

## V. CONCLUSION

In this paper, a novel bi-level energy management framework for pelagic island microgrid groups with joint energy and reserve market is proposed. In the developed framework, microgrid operators in the upper level aim to minimize the cost of the system and RSA are introduced to formulate the lower level. Due to the uneven temporal and spatial distribution of renewable energy resources and loads among different pelagic islands, RSA can practice arbitrage by paying service fees for vessels sailing among the islands to transmit batteries in the ESS. Chance constrained programming is applied in the framework to formulate the uncertainty and incorporate with the reserve service. Because of the integer variables such as numbers of vessels and batteries swapping and nonlinear terms in the lower level, a novel upper bound tightening reformulation and decomposition is proposed to solve the bi-level mixed-integer nonlinear problem for the operation frameworks. Case studies demonstrate the effectiveness of the proposed framework and show great edges on optimality gaps and computational efficiency by comparisons with existing algorithms. Future work will extend the proposed centralized approach to large-scale system with more energy exchange channels and investigate the improved centralized approach with a high computation efficiency. Furthermore, privacy-preserve approaches for the large-scale pelagic island microgrid groups will be explored with no private information exchange among different stakeholders.

## APPENDIX

**Lemma 1**

If the MP remains feasible, the algorithm generates integer variables combination $y_{l,k,i,t}$ to be added to $Z_k^l$ respectively at iteration $k$.

**Proof:**

From the reformulation process in (41)-(43), one can get that when $Z_k^l$ includes all the integer variables combination, the master problem in the iterative process is equivalent to the original bilevel problem. For any $Z_k^l$, the master problem is a relaxation of the original bilevel problem (See more detailed proofs in *Lemma 3.14* in [35]). In (46), (47), one can fix the integer variables as the parameters of the optimization problem in $Z_k^l$ and derive the primal and dual feasibility constraints for the continuous variables in the lower level through the strong duality theorem. Let $(x_{u,i,t}^*, x_{l,i,t}^*, y_{l,i,t}^*)$ be the solution of MP (47) at iteration $k$. As the SP (44) is feasible with the solution $x_{u,i,t}^*$, thus an integer variables combination is always generated in each iteration as long as the master problem remains feasible.

**Lemma 2**

If the generated integer variables combination $y_{l,k,i,t}$ already exists in the set $Z_k^l$ in iteration k, the proposed algorithm terminates with a bi-level optimal solution.

**Proof:**

If SP2 (45) is infeasible in iteration k, then $y_{l,k,i,t}$ is part of the optimal solution of SP computed in the previous iteration and consequently $\theta_k(x_{u,i,t}^*) = \theta_k(x_{u,i,t}^*, x_{l,i,t}, y_{l,i,t})$ holds. If SP2 (45) is feasible at iteration $k$, then $y_{l,k,i,t}$ is part of its optimal solution and also satisfies $\theta_k(x_{u,i,t}^*) = \theta_k(x_{u,i,t}^*, x_{l,i,t}, y_{l,i,t})$. If the generated integer variables combination already exists in the $Z_k^l$ and the corresponding constraints are also already contained in the MP (47) at iteration $k$. In particular, the constraint

$$f_{l,i,t}(x_{u,i,t}, x_{l,i,t}, y_{l,i,t}) \leq \min\{f_{l,i,t}(x_{u,i,t}, x_{c,l,i,t}, y_{c,l,i,t}) : Dx_{c,l,i,t} + Ey_{c,l,i,t} \leq \mathbf{b}_l\}$$

in (43), the primal problem of MP, is active at iteration k. As $\theta_k(x_{u,i,t}^*) = \theta_k(x_{u,i,t}^*, x_{l,i,t}, y_{l,i,t})$, the constraints for the upper-level decision variables $x_{u,i,t}^*$ corresponding to $y_{l,k,i,t}$ is exactly the optimal-value-function constraint, cf. Section 1.2 in [40]. Therefore, $(x_{u,i,t}^*, x_{l,i,t}^*, y_{l,i,t}^*)$ is a bi-level feasible solution. According to the *Lemma 3.14* in [35], it is also the optimal solution to the original bi-level problem.

In conclusion, the SP can generate integer variables combination $y_{l,k,i,t}$ and the corresponding constraints for MP at iteration $k$ in the decomposition process based on **Lemma 1** and MP can obtain the optimal solution when the generated integer variables combination $y_{l,k,i,t}$ from SP already exists in the set $Z_k^l$ at iteration $k$ according to **Lemma 2**.

REFERENCES

[1] Y. Han, K. Zhang, H. Li, E. A. A. Coelho and J. M. Guerrero, "MAS-Based distributed coordinated control and optimization in microgrid and microgrid clusters: a comprehensive overview," *IEEE Trans. Power Electronics*, vol. 33, no. 8, pp. 6488-6508, Aug. 2018.

[2] Z. Wang, K. Yang and X. Wang, "Privacy-preserving energy scheduling in microgrid systems," *IEEE Trans. Smart Grid*, vol. 4, no. 4, pp. 1810-1820, Dec. 2013.

[3] X. Sun, H. Xie, Y. Wang, C. Chen and Z. Bie, "Pricing for TSO-DSO coordination: a decentralized incentive compatible approach," *IEEE Trans. Power Systems*, doi: 10.1109/TPWRS.2022.3170436, 2022.

[4] M. Fathi and H. Bevrani, "Statistical cooperative power dispatching in interconnected microgrids," *IEEE Trans. Sustain. Energy*, vol. 4, no. 3, pp. 586-593, July 2013.

[5] R. Lahon, C. P. Gupta and E. Fernandez, "Optimal power scheduling of cooperative microgrids in electricity market environment," *IEEE Trans. Industrial Informatics*, vol. 15, no. 7, pp. 4152-4163, July 2019.

[6] H. Kim, J. Lee, S. Bahrami and V. W. S. Wong, "Direct energy trading of microgrids in distribution energy market," *IEEE Trans. Power Syst.*, vol. 35, no. 1, pp. 639-651, Jan. 2020.

[7] F. Pacaud, M. De Lara, J. -P. Chancelier and P. Carpentier, "Distributed multistage optimization of large-scale microgrids under stochasticity," *IEEE Trans. Power Syst.*, vol. 37, no. 1, pp. 204-211, Jan. 2022.

[8] W. Huang, W. Zheng and D. J. Hill, "Distributionally robust optimal power flow in multi-microgrids with decomposition and guaranteed convergence," *IEEE Trans. Smart Grid*, vol. 12, no. 1, pp. 43-55, Jan. 2021.

[9] D. Chen, K Chen, Z. Li, et al., "PowerNet: multi-agent deep reinforcement learning for scalable powergrid control," *IEEE Trans. Power Syst.*, vol. 37, no. 2, pp. 1007-1017, March 2022.

[10] L. Che, M. Shahidehpour, A. Alabdulwahab and Y. Al-Turki, "Hierarchical coordination of a community microgrid with AC and DC microgrids," *IEEE Trans. Smart Grid*, vol. 6, no. 6, pp. 3042-3051, Nov. 2015.

[11] D. Papadaskalopoulos, D. Pudjianto and G. Strbac, "Decentralized coordination of microgrids with flexible demand and energy storage," *IEEE Trans. Sustain. Energy*, vol. 5, no. 4, pp. 1406-1414, Oct. 2014.

[12] X. Lin, C. Chen, and X. Zhou, "Integrated energy supply system of pelagic clustering islands, " *Proc. CSEE*, vol. 37, no. 1, pp. 98–110, 2017.

[13] M. Hu, Y. Wang, J. Xiao, and X. Lin, "Multi-energy management with hierarchical distributed multi-scale strategy for pelagic islanded microgrid clusters, " *Energy*, vol. 185, pp. 910–921, Oct. 2019.

[14] Q. Sui, F. Wei, C. Wu, X. Lin and Z. Li, "Day-ahead energy management for pelagic island microgrid groups considering non-integer-hour energy transmission," *IEEE Trans. Smart Grid*, vol. 11, no. 6, pp. 5249-5259, Nov. 2020.

[15] F. D. Kanellos, G. J. Tsekouras and N. D. Hatziargyriou, "Optimal demand-side management and power generation scheduling in an all-electric ship," *IEEE Trans. Sustain. Energy*, vol. 5, no. 4, pp. 1166-1175, Oct. 2014.

[16] K. Hein, Y. Xu, G. Wilson and A. K. Gupta, "Coordinated optimal voyage planning and energy management of all-electric ship with hybrid energy storage system," *IEEE Trans. Power Systems*, vol. 36, no. 3, pp. 2355-2365, May 2021.

[17] A. A. Anderson and S. Suryanarayanan, "Review of energy management and planning of islanded microgrids," *CSEE Journal of Power and Energy Systems*, vol. 6, no. 2, pp. 329-343, June 2020.

[18] Y. Wang, J. Qiu, Y. Tao and J. Zhao, "Carbon-Oriented Operational Planning in Coupled Electricity and Emission Trading Markets," *IEEE Trans. Power Systems*, vol. 35, no. 4, pp. 3145-3157, July 2020.

[19] H. Wu, M. Shahidehpour, Z. Li and W. Tian, "Chance-constrained day-ahead scheduling in stochastic power system operation," *IEEE Trans. Power Systems*, vol. 29, no. 4, pp. 1583-1591, July 2014.

[20] G. Wu, Y. Xiang, J. Liu, X. Zhang and S. Tang, "Chance-constrained optimal dispatch of integrated electricity and natural gas systems considering medium and long-term electricity transactions," *CSEE Journal of Power and Energy Systems*, vol. 5, no. 3, pp. 315-323, Sept. 2019.

[21] J. Li, Z. Wu, S. Zhou, H. Fu and X. -P. Zhang, "Aggregator service for PV and battery energy storage systems of residential building," *CSEE Journal of Power and Energy Systems*, vol. 1, no. 4, pp. 3-11, Dec. 2015.

[22] X. Duan, Z. Hu and Y. Song, "Bidding strategies in energy and reserve markets for an aggregator of multiple EV fast charging stations with battery storage," *IEEE Trans. Intelligent Transportation Systems*, vol. 22, no. 1, pp. 471-482, Jan. 2021.

[23] Z. Yi, Y. Xu, J. Zhou, W. Wu and H. Sun, "Bi-level programming for optimal operation of an active distribution network with multiple virtual power plants," *IEEE Trans. Sustain. Energy*, vol. 11, no. 4, pp. 2855-2869, Oct. 2020.

[24] Y. Lu, X. Yu, X. Jin, H. Jia and Y. Mu, "Bi-level optimization framework for buildings to heating grid integration in integrated community energy systems," *IEEE Trans. Sustain. Energy*, vol. 12, no. 2, pp. 860-873, April 2021.

[25] L. Yu, T. Zhao, Q. Chen and J. Zhang, "Centralized bi-level spatial-temporal coordination charging strategy for area electric vehicles," *CSEE Journal of Power and Energy Systems*, vol. 1, no. 4, pp. 74-83, Dec. 2015.

[26] A. Paudel, K. Chaudhari, C. Long and H. B. Gooi, "Peer-to-peer energy trading in a prosumer-based community microgrid: a game-theoretic model," *IEEE Trans. Ind. Electron.*, vol. 66, no. 8, pp. 6087-6097, Aug. 2019.

[27] T. Zhao, J. Zhang and P. Wang, "Flexible active distribution system management considering interaction with transmission networks using information-gap decision theory," *CSEE Journal of Power and Energy Systems*, vol. 2, no. 4, pp. 76-86, Dec. 2016.

[28] H. Liu, W. Wu and Y. Wang, "Bi-level off-policy reinforcement learning for two-timescale Volt/VAR control in active distribution networks," *IEEE Trans. Power Syst.*, doi: 10.1109/TPWRS.2022.3168700, 2022.

[29] Moore, James T., and Jonathan F. Bard. "The mixed integer linear bilevel programming problem. " *Operations Research*, vol. 38, no. 5, pp. 911–21, 1990.

[30] Domínguez, L. F. and E. N. Pistikopoulos. "Multiparametric programming based algorithms for pure integer and mixed-integer bilevel programming problems." *Computers & Chemical Engineering*, vol. 34, no. 5, pp. 2097-2106, 2010.

[31] Mitsos, A. "Global solution of nonlinear mixed-integer bilevel programs. " *J Glob Optim* vol. 47, pp. 557–582, 2010.

[32] Kleniati, P.-M. and C. S. Adjiman, "A generalization of the branch-and-sandwich algorithm: from continuous to mixed-integer nonlinear bilevel problems." *Computers & Chemical Engineering*, vol. 72, pp. 373-386, 2015.

[33] B. Zeng and Y. An, "Solving bilevel mixed integer program by reformulations and decomposition, " *Optimization online*, pp. 1–34, 2014.

[34] D. Yue, J. Gao, B. Zeng and F. You. "A projection-based reformulation and decomposition algorithm for global optimization of a class of mixed integer bilevel linear programs," *Journal of Global Optimization*, vol. 73(1), pp 27-57, January, 2018.

[35] Merkert, M., Orlinskaya, G. and Weninger, D. "An exact projection-based algorithm for bilevel mixed-integer problems with nonlinearities." *J Glob Optim*, vol. 84, pp. 607-650, 2022.

[36] Shen, Z., et al. "An accelerated Stackelberg game approach for distributed energy resource aggregator participating in energy and reserve markets considering security check." *International Journal of Electrical Power & Energy Systems*, vol. 142: 108376, 2022.

[37] Y. Ye, D. Papadaskalopoulos, J. Kazempour and G. Strbac, "Incorporating non-convex operating characteristics into bi-level optimization electricity market models," *IEEE Trans. Power Syst.*, vol. 35, no. 1, pp. 163-176, Jan. 2020.

[38] C. Ruiz, A. J. Conejo, and S. A. Gabriel, "Pricing non-convexities in an electricity pool," *IEEE Trans. Power Syst.*, vol. 27, no. 3, pp. 1334–1342, Aug. 2012.

[39] D. Huppmann and S. Siddiqui, "An exact solution method for binary equilibrium problems with compensation and the power market uplift problem," Eur . J. Oper . Res., vol. 266, no. 2, pp. 622–638, Apr. 2018.

[40] Dempe, S., Kalashnikov, V., Prez-Valds, G.A., Kalashnykova, N. "Bilevel Programming Problems: Theory, Algorithms and Applications to Energy Networks," Springer, Berlin, 2015.